\documentclass[a4paper,12pt]{article}

\usepackage[top=3.0cm,bottom=3.0cm,left=2.45cm,right=2.45cm]{geometry}
\usepackage{amsmath,amsthm,mathrsfs,graphicx,amsfonts}
\usepackage{bm}
\usepackage{cases}
\usepackage{hyperref}
\usepackage[english]{babel}
\usepackage{amsmath,amsthm}
\usepackage{amsfonts}
\usepackage{latexsym}
\usepackage{graphicx}
\usepackage{txfonts}
\usepackage[numbers,sort&compress]{natbib}
\usepackage[natural]{xcolor}
\usepackage{rotating}
\usepackage{mathtools}
\usepackage{enumitem}
\usepackage{appendix}
\allowdisplaybreaks
\newtheorem{theorem}{\bf Theorem}[section]
\newtheorem{lemma}[theorem]{Lemma}

\newtheorem{corollary}[theorem]{Corollary}

\title{Tight bounds towards Zarankiewicz problem in hypergraph}
\author{Guorong Gao$^{a,b}$\thanks{
Research supported by National Key R\&D Program of China (Grant No. 2023YFA1010202),
National Natural Science Youth Foundation of China (Grant No. 12401448), Natural Science Foundation of Fujian Province (Grant No. 2024J08030).
Email: grgao@fzu.edu.cn}, Jianfeng Hou$^{b}$\thanks{
Research supported by National Key R\&D Program of China (Grant No. 2023YFA1010202), National Natural Science Foundation of China (Grant No. 12071077), the Central Guidance on Local Science and Technology Development Fund of Fujian Province (Grant No. 2023L3003).
Email: jfhou@fzu.edu.cn},\, Shuping Huang$^c$\thanks{Email: hsp@mail.ustc.edu.cn}, Hezhi Wang$^{a}$\thanks {
Email: sdlgsjwhz@126.com }\\
\small$^a$Center for Discrete Mathematics and Theoretical Computer Science,\\
\small Fuzhou University, Fuzhou, Fujian, China\\
\small$^b$School of Mathematics and Statistics,
\small Fuzhou University, Fuzhou, Fujian, China\\
\small$^c$School of Mathematical Sciences, University of Sciences and Technology of China, Hefei, China}
\date{}

\begin{document}

\maketitle

\date{}
\maketitle

\begin{abstract}
The classical Zarankiewicz problem, which concerns the maximum number of edges in a bipartite graph without a forbidden complete bipartite subgraph, motivates a direct analogue for hypergraphs. Let $K_{s_1,\ldots, s_r}$ be the complete $r$-partite $r$-graph such that the $i$-th part has $s_i$ vertices. We say an $r$-partite $r$-graph $H=H(V_1,\ldots,V_r)$ contains an ordered $K_{s_1,\ldots, s_r}$ if $K_{s_1,\ldots, s_r}$ is a subgraph of $H$ and the set of size $s_i$ vertices is embedded in $V_i$.
The Zarankiewicz number for $r$-graph, denoted by $z(m_1, \ldots, m_{r}; s_1,, \ldots,s_{r})$, is the maximum number of edges of the $r$-partite $r$-graph whose $i$-th part has $m_i$ vertices and does not contain an ordered $K_{s_1,\ldots, s_r}$. In this paper, we show that
$$z(m_1,m_2, \cdots, m_{r-1},n ; s_1,s_2, \cdots,s_{r-1}, t)=\Theta\left(m_1m_2\cdots m_{r-1} n^{1-1 / s_1s_2\cdots s_{r-1}}\right)$$
for a range of parameters. This extends a result of Conlon [Math. Proc. Camb. Philos. Soc. (2022)].
\end{abstract}


\section{Introduction}

An $r$-uniform hypergraph (or $r$-graph for convenience) $H = (V(H),E(H))$ is a pair of a vertex set $V(H)$ and an edge set $E(H)$, where the edge set is a collection of $r$-element subsets of the vertex set.

Given a graph $F$, the Tur\'an number of $F$, denoted by $\operatorname{ex}(n, F)$, is the maximum possible number of edges in the $n$-vertex $F$-free graphs. The classical Erd\H{o}s-Stone-Simonovits theorem  gives an estimate for this function, showing that
$$
\operatorname{ex}(n, F)=\left(1-\frac{1}{\chi(F)-1}+o(1)\right)\binom{n}{2},
$$
where $\chi(F)$ is the chromatic number of $F$. For bipartite $F$, this gives the bound $\operatorname{ex}(n, F)=o\left(n^2\right)$. While more precise estimates are known, a number of notoriously difficult open problems remain.
The most intensively studied case is when $F=K_{s, t}$, the complete bipartite graph with parts of order $s$ and $t$.
In the 1930s, Erd\H{o}s--R\'enyi--Brown \cite{B66} gave optimal lower bounds for $K(2,t)$ and $K(3,t)$. Then Koll\'ar--R\'onyai--Szab\'o \cite{KRT96} and Alon--R\'onyai--Szab\'o \cite{ARS99} proved that
\begin{align}
\operatorname{ex}\left(n, K_{s,t}\right)=\Theta\left(n^{r-\frac{1}{s}}\right)
\end{align}
where $t>(s-1)!$. Bukh \cite{B15,B24} introduced the powerful random algebraic method and showed that (1) holds as long as $t>9^{s+o(s)}$.

Write $G= G(m, n)$  with parts of size $m$ and $n$. The Zarankiewicz number $z(m, n; s, t)$ is the maximum number of edges in $G(m, n)$ with parts $U$ and $V$ satisfying that there is no copy of $K_{s,t}$ with $s$ vertices in $U$ and $t$ vertices in $V$. The classical K\H{o}v\'ari-S\'os-Tur\'an Theorem \cite{KST54} gives
$$z(m, n ; s, t)=O\left(m n^{1-\frac{1}{s}}\right).$$
By the random algebraic method, Conlon \cite{C22} proved that
$$z(m, n ; s, t)=\Theta\left(m n^{1-\frac{1}{s}}\right)$$
for any fixed $2\leq s \leq t$ and any $m \leqslant n^{t^{1 /(s-1)} / s(s-1)}$.

Let $K_{s_1,\ldots, s_r}$ be the complete $r$-partite $r$-graph such that the $i$-th part has $s_i$ vertices. Ma, Yuan and Zhang \cite{MYZ18} proved that $$\operatorname{ex}\left(n, K_{s_1, s_2, \ldots, s_{r-1}, t}\right)=\Omega\left(n^{r-\frac{1}{s_1 s_2 \cdots s_{r-1}}}\right)$$
for sufficiently large $t$. More recently, Pohoata and Zakharov \cite{PZ21} proved the same lower bound as long as $t>((r-1)(s-1))!$ and Mubayi \cite{M25} improved this lower bound on $t$ substantially in the Zarankiewicz case, from factorial to exponential at the expense of a small $o(1)$   error parameter in the exponent.

In this paper, we study the Zarankiewicz problem in hypergraph. Define that an $r$-partite $r$-graph $H=H(V_1,\ldots,V_r)$ contains an ordered $K_{s_1,\ldots, s_r}$ if $K_{s_1,\ldots, s_r}$ is a subgraph of $H$ and the set of size $s_i$ vertices is embedded in $V_i$.
The Zarankiewicz number for $r$-graph, denoted by $z(m_1, \ldots, m_{r}; s_1,, \ldots,s_{r})$, is the maximum number of edges of the $r$-partite $r$-graph whose $i$-th part has $m_i$ vertices and does not contain an ordered $K_{s_1,\ldots, s_r}$. If $m_1=m_2=\cdots=m_r=n$, then we write $z(n; {s_1, \ldots, s_{r}})$ for short. Recently, Mubayi \cite{M25} proved that
$$z \left(n; {s_1, \ldots, s_{r-1}, t}\right)>n^{1-o(1)} \cdot z \left(n; {s_1, \ldots, s_{r-3}, s_{r-2} s_{r-1}, t}\right)$$
for any $r \geq 3$, and positive integers $s_1, \ldots, s_{r-1}, t$, $n \to \infty$.

Our first result is a supersaturation result for $K_{s_1,\ldots, s_r}$, which extends a classical Theorem of Erd\H{o}s \cite{E64}.

\begin{theorem}\label{thmupper}
Let  $H$ be an $r$-partite $r$-graph with parts $V_1,\ldots V_{r}$, $\lvert V_i \rvert = m_i$ for $1\leq i\leq r$ and $m_r=\min_{1\leq i\leq r}\{m_i\}$. Then there exists constants $c_1$ and $c_2$, such that if $|E(H)| \geq c_1 \cdot\left(\prod_{i=1}^{r-1} m_i\right)\cdot m_r^{1-\frac{1}{s_1s_2\cdots s_{r-1}}}$, then $H$ contains at least $ c_2\cdot \prod_{i=1}^{r}\binom{ m_i}{s_i} \cdot p^{\prod_{i=1}^{r}s_i}$  copies of ordered $K_{s_1,s_2 \cdots s_{r}}$,
where $p=\frac{|E(H)|}{\prod_{i=1}^{r} m_i}$.
\end{theorem}

As a corollary of Theorem \ref{thmupper}, we have the following upper bound for the Zarankiewicz problem, which extends the K\H{o}v\'ari-S\'os-Tur\'an Theorem \cite{KST54}.
\begin{corollary}\label{coroupper}
     For any fixed $m_1,m_2,\cdots, m_{r}$ and $s_1,s_2, \cdots ,s_{r}$, if $m_i\geq m_r$ for all $1\leq i\leq r-1$, then
	$$
	z(m_1,m_2, \cdots, m_{r} ; s_1,s_2, \cdots,s_{r})=O\left(m_1m_2\cdots m_{r-1} m_r^{1-1 / s_1s_2\cdots s_{r-1}}\right) .
	$$
\end{corollary}

By a generalization of Bukh's random algebraic method, we prove the following lower bound for the Zarankiewicz problem, which extends the result of Conlon \cite{C22}.

\begin{theorem}\label{thmlower}
	
For any fixed $m_1,m_2,\cdots,m_{r-1}, n$ and $s_1,s_2, \cdots ,s_{r-1},t$, let $m=\prod_{i=1}^{r-1} m_i$ and $s=\prod_{i=1}^{r-1} s_i$. If $s \leq t$ and $m  \leq n^{t^{1 /s-1} / s (s-1)}$, then
$$
z(m_1,m_2, \cdots, m_{r-1},n ; s_1,s_2, \cdots,s_{r-1}, t)=\Omega\left(m_1m_2\cdots m_{r-1} n^{1-1 / s_1s_2\cdots s_{r-1}}\right) .
$$
\end{theorem}

Combining Corollary \ref{coroupper} and Theorem \ref{thmupper}, we obtain a range of tight bound (up to the constant) for Zarankiewicz problem in hypergraph.

\begin{corollary}
For any fixed $m_1,m_2,\cdots,m_{r-1}, n$ and $s_1,s_2, \cdots ,s_{r-1},t$, let $m=\prod_{i=1}^{r-1} m_i$ and $s=\prod_{i=1}^{r-1} s_i$.
If $s \leq t$, $m  \leq n^{t^{1 /s-1} / s (s-1)}$  and $m_i\geq n$ for all $1\leq i\leq r-1$, then
	$$
	z(m_1,m_2, \cdots, m_{r-1},n ; s_1,s_2, \cdots,s_{r-1}, t)=\Theta\left(m_1m_2\cdots m_{r-1} n^{1-1 / s_1s_2\cdots s_{r-1}}\right) .
	$$
\end{corollary}

In Section 2, we prove the Theorem \ref{thmlower}. Theorem \ref{thmupper} is proved in Section 3. Throughout this paper, we adopt the standard notations. In particular, for an $r$-graph $H$, the degree of a vertex $v\in V(H)$, denoted by $d_H(v)$, is the number of edges of $H$ containing $v$. We omit the subscript if there is no confusion.
\vskip 0.5cm
\noindent{\textbf{Remark.}} After we completed this work, we noticed that Chen, Liu and Ye \cite{CLY25}  established the same lower bound of $z(m_1,m_2, \cdots, m_{r-1},n ; s_1,s_2, \cdots,s_{r-1}, t)$ for a larger $t$ and without restriction of $n$. Their method is an another generalization of Bukh's random algebraic method.

\section{Proof of the Theorem \ref{thmlower}}

\hspace{0.5cm}
The proof employs the powerful random algebraic method, which is developed in \cite{B15,B24,C22}.
For a prime power $q$, let $\mathbb{F}_q$ be the finite field of order $q$. Writing a polynomial in $t$ variables over $\mathbb{F}_q$ as $f(X)$: ${\mathbb{F}}^{t}_q \rightarrow {\mathbb{F}}_q$, where $X=\left(X_1, \ldots, X_t\right)$.

Let $\mathcal{P}_d$ be the set of polynomials in $X$ of degree at most $d$, that is, the expression for $f$ is as follows:
$$
\sum_{\sum_{i=1}^ta_i \leq d} \alpha X_1^{a_1} X_2^{a_2} \cdots X_t^{a_t} \quad \alpha  \in  {\mathbb{F}}_q.
$$
  By a random polynomial, we just mean a polynomial chosen uniformly from the set $\mathcal{P}_d$. One may produce such a random polynomial by taking the coefficients of the monomials above to be independent random elements of $\mathbb{F}_q$.

The following lemma, due to Conlon \cite{C22}, estimates the probability that a randomly selected polynomial from $\mathcal{P}_d$ passes through $m$ distinct points. Let $\overline{\mathbb{F}}_q$ be the algebraic closure of $\mathbb{F}_q$.

\begin{lemma}[\cite{C22}]\label{lem1}
	Suppose that $q>\binom{m}{2}$ and $d \geq m-1$. If $f$ is a random $t$-variate polynomial of degree $d$ over $\mathbb{F}_q$ and $x_1, \ldots, x_m$ are $m$ distinct points in $\overline{\mathbb{F}}_q^t$, then
	$$
	\mathbb{P}\left[f\left(x_i\right)=0 \text { for all } i=1, \ldots, m\right] \leq 1 / q^m .
	$$
\end{lemma}

Over an algebraically closed field $\overline{\mathbb{F}}$, a variety is a set of the form

$$
W=\left\{x \in \overline{\mathbb{F}}^t: f_1(x)=\cdots=f_s(x)=0\right\}
$$
for some collection of polynomials $f_1, \ldots, f_s: \overline{\mathbb{F}}^t \rightarrow \overline{\mathbb{F}}$. The variety is irreducible if it cannot be written as the union of two proper subvarieties. The dimension $\operatorname{dim} W$ of $W$ is then
the maximum integer $d$ such that there exists a chain of irreducible subvarieties of $W$ of the form
$$
\emptyset \subsetneq\{p\} \subsetneq W_1 \subsetneq W_2 \subsetneq \cdots \subsetneq W_d \subset W,
$$
where $p$ is a point.  Another definition of the dimension of an algebraic variety is that
$$
\operatorname{dim} W=\max \left\{\operatorname{dim} W_i \mid W_i \ is \ an \ irreducible \ component \ of \ W       \right\}.
$$

In what follows we introduce three standard lemmas about the varieties.
\begin{lemma}\label{lem2}
	Every variety $W$ over an algebraically closed field $\overline{\mathbb{F}}$ with $\operatorname{dim} W \geq 1$ has infinitely many points.
\end{lemma}
\begin{lemma}\label{lem3}
	Suppose that $W$ is an irreducible variety over an algebraically closed field $\overline{\mathbb{F}}$. Then, for any polynomial $g: \overline{\mathbb{F}}^t \rightarrow \overline{\mathbb{F}}, W \subseteq\{x: g(x)=0\}$ or $W \cap\{x: g(x)=0\}$ is a variety of dimension less than $\operatorname{dim} W$.
\end{lemma}
\begin{lemma}[B\'ezout's theorem \cite{F84}] \label{lem4}
If, for a collection of polynomials $f_1, \ldots, f_t$ : $\overline{\mathbb{F}}^t \rightarrow \overline{\mathbb{F}}$, the variety
	$$
	W=\left\{x \in \overline{\mathbb{F}}^t: f_1(x)=\cdots=f_t(x)=0\right\}
	$$
	has $\operatorname{dim} W=0$, then
	$$
	|W| \leq \prod_{i=1}^t \operatorname{deg}\left(f_i\right).
	$$
	Moreover, for a collection of polynomials $f_1, \ldots, f_s: \overline{\mathbb{F}}^t \rightarrow \overline{\mathbb{F}}$, the variety
	$$
	W=\left\{x \in \overline{\mathbb{F}}^t: f_1(x)=\cdots=f_s(x)=0\right\}
	$$
	has at most $\prod_{i=1}^s \operatorname{deg}\left(f_i\right)$ irreducible components.
\end{lemma}

\begin{proof}[Proof of Theorem \ref{thmlower}]
Let $l=\prod_{i=1}^{r-1} l_i$, $s=\prod_{i=1}^{r-1} s_i$, fix $d=\left\lceil t^{1 /(s-1)}\right\rceil-1$ and $\ell=\left\lfloor q^{(d+1) /(s-1)} / 2 d\right\rfloor$. Consider the $r$-partite $r$-graph between sets $U_1,U_2  \cdots U_{r-1}$ and $V$, where $V$ may be viewed as a copy of $\mathbb{F}_q^{s}$ for some prime power $q$ and each $U_i$ has order $\ell_i$, $i=1,2 \cdots r-1$, every vertex class
$$
\left(u_{p ^1}, u_{p ^2}, \cdots, u_{p ^{r-1}}\right), \quad u_{p ^1} \in U_1, u_{p ^2} \in U_2 \cdots u_{p ^{r-1}} \in U_{r-1},
$$
of which is associated to an $\left(\prod_{i=1}^{r-1} s_i-1\right)$-variate polynomial $f_{p^1, \cdots, p^{r-1}}$ of degree at most $d$ with coefficients in $\mathbb{F}_q$. Each $\left(u_{p ^1}, u_{p ^2}, \cdots, u_{p ^{r-1}}\right)$ is then joined to the set of points
$$
S_{p^1, \cdots, p^{r-1}}=\left\{\left(x_1, \ldots, x_{\prod_{i=1}^{r-1} s_i-1}, f_{p^1, \cdots, p^{r-1}}\left(x_1, \ldots, x_{\prod_{i=1}^{r-1} s_i-1}\right)\right): x_1, \ldots, x_{\prod_{i=1}^{r-1} s_i-1} \in \mathbb{F}_q\right\}
$$
in $V$. Note that, for any $1 \leq j_i \leq s_i$ and $j=\prod_{i=1}^{r-1}j_i$  ($1 \leqslant p_{j_i}^i \leqslant l_i$) elements of
$$
\left(p^1, p^2, \cdots, p^{r-1}\right), \quad 1 \leq p^1 \leq  l_1, \cdots, 1 \leq p^{r-1} \leq l_{r-1},
$$
we denote them by $a_1, a_2, \cdots, a_j$.
$$
S_{a_1} \cap \cdots \cap S_{a_j}=\left\{\left(x_1, \ldots, x_s\right): x_s=f_{a_1}\left(x_1, \ldots, x_{\prod_{i=1}^{r-1} s_i-1}\right)=\cdots=f_{a_j}\left(x_1, \ldots, x_{\prod_{i=1}^{r-1} s_i-1}\right)\right\}
$$
This intersection therefore has the same size as $T_{a_1, a_2} \cap \cdots \cap T_{a_1, a_j}$, where
$$
T_{a_i, a_i^{\prime}}=\left\{\left(x_1, \ldots, x_{\prod_{i=1}^{r-1} s_i-1}\right):\left(f_{a_i}-f_{a^{\prime}_i}\right)\left(x_1, \ldots, x_{\prod_{i=1}^{r-1} s_i-1}\right)=0\right\} .
$$

\begin{lemma}\label{lem5}
Let $H$ is a $r$-partite $r$-graph   with parts $V_1,V_2$,  $\cdots$, $V_{r-1} $ and $U$, $\lvert V_i \rvert = l_i , \lvert U \rvert =n$, $i=1,2, \cdots, r-1$, let $l=\prod_{i=1}^{r-1} l_i$, $s=\prod_{i=1}^{r-1} s_i$, $j=\prod_{i=1}^{r-1} j_i$, there exist a choice of  $  f_{p^1, \cdots, p^{r-1}}$ of degree at most $d$ with $1 \leq p^1 \leq  l_1, \cdots, 1 \leq p^{r-1} \leq l_{r-1}$ such that the dimension of the intersection $T_{a_1, a_2} \cap \cdots \cap T_{a_1, a_j}$ is at most $s-j$.
\end{lemma}

The main idea of the proof of Lemma \ref{lem5} can be found in \cite{C22}. However, it is not a separate proof. To ensure the completeness of the proof, we include the proof of this lemma in the appendix.

By Lemma \ref{lem5}, there is a choice of $f_{p^1, \cdots, p^{r-1}}$ for  $1 \leq p^1 \leq  l_1, \cdots, 1 \leq p^{r-1} \leq l_{r-1}$ such that  for any $1 \leq j_i \leq s_i$ and  $1\leq p^i_1<p^i_2< \cdots <p^i_{j_i} \leq l_i$, $i=1,2, \cdots ,r-1$, the intersection $T_{a_1, a_2} \cap \cdots \cap T_{a_1, a_s}$ has dimension  $0$, so by Lemma 2.4, the number of points in the intersection is at most $d^{s-1}<t $. Therefore, for any   $1 \leq p^1 \leq  l_1, \cdots, 1 \leq p^{r-1} \leq l_{r-1}$, the intersection $S_{a_1} \cap \cdots \cap S_{a_s}$ has at most $t-1$ points, so there is no copy of $K_{s_1,s_2, \cdots, s_{r-1},t}$ with $s_i$ vertices in $U_i$ and $t$ vertices in $V$. Since $l=\Omega\left(q^{\frac{d+1}{s-1}}\right)$, $|V|=q^s$ and $|E|=l q^{s-1}$, so for $m_0:=n^{\frac{d+1}{s(s-1)}} \leqslant n^{t^{1 /s-1} / s (s-1)}$, $m_1m_2\cdots m_{r-1}\leqslant m_0$, we have the following result,
$$
	z(m_1,m_2, \cdots, m_{r-1},n ; s_1,s_2, \cdots, s_{r-1}, t)=\Omega\left(m_1m_2\cdots m_{r-1} n^{1-1 / s_1s_2\cdots s_{r-1}}\right) .
$$
By applying Bertrand's postulate, we can extend this result to all $n$ not only the form $q^s$ with $q$ a prime power.
\end{proof}

\section{Proof of the Theorem \ref{thmupper}}
\begin{proof}[Proof of Theorem \ref{thmupper}]
We will prove it by induction on $r$.
Set $c_1$ be a sufficiently large constant and $c_2$ be a sufficiently small constant.
For a number $x\geq 0$ and a positive integer s, let
\[\binom{x}{s}=\left\{\begin{array}{cl}{0,}&{\quad\mathrm{if}~x<s-1;}\\{\frac{x(x-1)\cdots(x-s+1)}{s!},}&{\quad\mathrm{if}~x\geq s-1.}\\\end{array}\right.\]

For the base case $r=2$, let $H$ be a bipartite graph with parts $V_1,V_2$ and $| V_1 | = m_1 , | V_2 | = m_2$. We use the standard double counting to prove that if $e=|E(H)| \geq c_1m_1 m_2^{1-\frac{1}{s_1}}$, then there are many copies of ordered $K_{s_1, s_2}$.

Let $t_{s_1, 1}$ be the number of ordered $K_{s_1, 1}$ of $H$ and  $t_{s_1, s_2}$ be the number of ordered $K_{s_1, s_2}$. Clearly,
\begin{equation*}
	t_{s_1, 1}=\sum_{v \in V_2}\binom{d(v)}{s_1}\geq m_2\binom{\frac{\sum_{v \in V_2}d(v)}{m_2}}{s_1}=m_2\binom{\frac{e}{m_2}}{s_1},
\end{equation*}
where the inequality follows from the Jensen's inequality.
Let $S\subseteq V_1$ be a vertex set of size $s_1$ and
$f(S)$ be the number of vertices adjacent to all vertices of $S$. Then
$\sum_{S \subset V_1 } f(S)=t_{s_1, 1}$ and we have

\begin{align*}
t_{s_1, s_2}=\sum_{S \subset V_1  }\binom{f(S)}{s_2}
&\geq \binom{ m_1}{s_1}\binom{{\sum_{S \subset V_1    } f(S)}/{\binom{ m_1}{s_1}}}{s_2} \\
&=\binom{m_1}{s_1}\binom{{t_{s_1, 1}}/{\binom{m_1}{s_1}}}{s_2}  \\
&\geq \binom{m_1}{s_1}\binom{{m_2\binom{\frac{e}{m_2}}{s_1}}/{\binom{m_1}{s_1}}}{s_2}.
\end{align*}
As $e \geq c_1m_1 m_2^{1-\frac{1}{s_1}}$ and $c_1$ is sufficiently large, we have ${m_2\binom{\frac{e}{m_2}}{s_1}}/{\binom{m_1}{s_1}}> s_2-1$. Thus
$$
t_{s_1, s_2}\geq c_2\binom{m_1}{s_1}\binom{m_2}{s_2}\left(\frac{e}{m_1 m_2}\right)^{s_1 s_2},
$$
where $c_2$ is a sufficiently small constant.

Suppose the theorem holds for $r-1$, we are going to prove the theorem holds for $r$.
Let $H$ be an $r$-partite $r$-graph with parts $V_1$, $V_2$, $\cdots$ ,$V_r$ and $|V_i|=m_i$ for $1\leq i\leq r$. To prove the theorem, we set
$m_i\geq m_r$ for all $1\leq i\leq r-1$ and $e=|E(H)|\geq c_1 \cdot\left(\prod_{i=1}^{r-1} m_i\right)\cdot m_r^{1-\frac{1}{s_1s_2\cdots s_{r-1}}}$.

Without loss of generality, let $m_{r-1}=\min_{1\leq i\leq r-1}\{m_i\}$. Remove all the vertices of $V_r$ with degree less than $c_1 \cdot\left(\prod_{i=1}^{r-2} m_i\right)\cdot m_{r-1}^{1-\frac{1}{s_1s_2\cdots s_{r-2}}}$. This process removes at most
$$m_r\times c_1 \cdot\left(\prod_{i=1}^{r-2} m_i\right)\cdot m_{r-1}^{1-\frac{1}{s_1s_2\cdots s_{r-2}}}\leq c_1 \cdot\left(\prod_{i=1}^{r-1} m_i\right)\cdot m_r^{1-\frac{1}{s_1s_2\cdots s_{r-2}}}=o(e)$$
edges, where the inequality follows from $ m_{r-1}\geq m_r$. Let $V_r'$ be the collection of the remaining vertices of $V_r$.
Then
$$\sum_{v\in V_r'} d(v)=(1-o(1))e\geq \frac{e}{2}.$$

For each $v\in V_r'$, let $H_v$ be link-hypergraph of $v$ with edge set
$\{h\backslash\{v\}:~v\in h\in E(H)\}$. Clearly, $|E(H_v)|=d(v)\geq c_1 \cdot\left(\prod_{i=1}^{r-2} m_i\right)\cdot m_{r-1}^{1-\frac{1}{s_1s_2\cdots s_{r-2}}}$ and $H_v$ is an ($r-1$)-partite ($r-1$)-graph. By induction, $H_v$ has at least $ c_2{'}\cdot \prod_{i=1}^{r-1}\binom{ m_i}{s_i} \cdot \left(\frac{d(v)}{\prod_{i=1}^{r-1} m_i}\right)^{\prod_{i=1}^{r-1}s_i}$ copies of ordered $K_{s_1,s_2, \cdots, s_{r-1}}$. Let $t_a$ be the number of ordered $K_{s_1,s_2, \cdots, s_{r-1},1}$ of $H$ and $t_b$ be the number of ordered $K_{s_1,s_2, \cdots, s_{r-1},s_r}$. Then by Jensen's inequality
\begin{align*}
t_a &\geq  \sum_{v \in V_r'}  c_2{'}\cdot \prod_{i=1}^{r-1}\binom{ m_i}{s_i} \cdot \left(\frac{d(v)}{\prod_{i=1}^{r-1} m_i}\right)^{\prod_{i=1}^{r-1}s_i}\\
&\geq c_2'\prod_{i=1}^{r-1}\binom{ m_i}{s_i} \cdot |V_r'|\cdot \left(\frac{\sum_{v\in V_r'}d(v)}{|V_r'|\cdot\prod_{i=1}^{r-1} m_i}\right)^{\prod_{i=1}^{r-1}s_i}\\
&\geq c_2'\prod_{i=1}^{r-1}\binom{ m_i}{s_i} \cdot m_r\cdot \left(\frac{e/2}{\prod_{i=1}^{r} m_i}\right)^{\prod_{i=1}^{r-1}s_i}.
\end{align*}

Let $\vec{S}$=$(S_1,S_2, \cdots ,S_{r-1})$, where $S_i \subset V_i$ and $|S_i|=s_i$ for $1\leq i\leq r-1$, we denote by $f(\vec{S})$ the number of vertices $v$ such that $ S_1 \cup S_2 \cup \cdots \cup S_{r-1}\cup v$ induces a copy of  ordered $K_{s_1,s_2, \cdots, s_{r-1},1}$ in $H$. Then $\sum_{\vec{S}} f(\vec{S})=t_a$. Again by Jensen's inequality,

\begin{align*}
t_b=\sum_{\vec{S}}\binom{f(\vec{S})}{s_r}
&\geq \prod_{i=1}^{r-1}\binom{ m_i}{s_i}\binom{{\sum_{\vec{S}} f(\vec{S})}/{\prod_{i=1}^{r-1}\binom{ m_i}{s_i}}}{s_r} \\
& = \prod_{i=1}^{r-1}\binom{ m_i}{s_i}\binom{{t_a}/{\prod_{i=1}^{r-1}\binom{ m_i}{s_i}}}{s_r} \\
&\geq \prod_{i=1}^{r-1}\binom{ m_i}{s_i}\binom{{ c_2'\prod_{i=1}^{r-1}\binom{ m_i}{s_i} \cdot m_r\cdot \left(\frac{e/2}{\prod_{i=1}^{r} m_i}\right)^{\prod_{i=1}^{r-1}s_i}}/{\prod_{i=1}^{r-1}\binom{ m_i}{s_i}}}{s_r} \\
&= \prod_{i=1}^{r-1}\binom{ m_i}{s_i}\binom{{ c_2' m_r\cdot \left(\frac{e/2}{\prod_{i=1}^{r} m_i}\right)^{\prod_{i=1}^{r-1}s_i}}}{s_r}
\end{align*}
As $e\geq c_1 \cdot\left(\prod_{i=1}^{r-1} m_i\right)\cdot m_r^{1-\frac{1}{s_1s_2\cdots s_{r-1}}}$ and $c_1$ is sufficiently large, we have ${ c_2' m_r\cdot \left(\frac{e/2}{\prod_{i=1}^{r} m_i}\right)^{\prod_{i=1}^{r-1}s_i}}>s_r-1$. Thus
$$ t_b\geq  c_2 \prod_{i=1}^{r}\binom{ m_i}{s_i} \left(\frac{e}{\prod_{i=1}^{r}m_i}\right)^{s_1s_2\cdots s_r},$$
as desired.

\end{proof}

\appendix

\section*{Appendix A. Proof of Lemma 2.5}
\begin{proof}[Proof of Lemma 2.5]
We are going to show that there is a choice of $f_{p^1, \cdots, p^{r-1}}$ for  $1 \leq p^1 \leq  l_1, \cdots, 1 \leq p^{r-1} \leq l_{r-1}$ such that  for any $1 \leq j_i \leq s_i$, $1 \leq k_i \leq l_i$ and  $1\leq p^i_1<p^i_2< \cdots <p^i_{j_i} \leq k_i \leq l_i$, $i=1,2, \cdots ,r-1$, the intersection $T_{a_1, a_2} \cap \cdots \cap T_{a_1, a_j}$ has dimension at most  $s-j$. To do this, we will pick the $f_{p^1, \cdots, p^{r-1}} $ in sequence and show by induction. For convenience, we let $f_{a_{1}}, f_{a_{2}},\cdots, f_{a_{l}}$ denote $l$ elements of the form for $f_{p^1, \cdots, p^{r-1}}$, $1 \leq p^1 \leq  l_1, \cdots, 1 \leq p^{r-1} \leq l_{r-1}$.

To begin the induction, we let $f_{a_{1}}$ be any $(s_1s_2 \cdots s_{r-1}-1)$-variate polynomial of degree $d$. In this case, the condition that the intersection $T_{a_1, a_2} \cap \cdots \cap T_{a_1, a_j} $ have dimension at most at most $s-j$, for any $1 \leq j_i \leq s_i$ and  $1\leq p^i_1<p^i_2< \cdots <p^i_{j_i} \leq k_i \leq l_i$, $i=1,2, \cdots ,r-1$ is degenerate, but can be meaningfully replaced by the observation that the set of all $\left(x_1, \ldots, x_{s-1}\right)$, corresponding to the trivial intersection, equals $\overline{\mathbb{F}}_q^{s-1}$, which has dimension $s_1s_2 \cdots s_{r-1}-1$, as require.

Let $k=\prod_{i=1}^{r-1}k_i$, suppose that $f_{a_1}, \ldots, f_{a_{k-1}}$ have been chosen consistent with the induction hypothesis. We would like to pick $f_{a_k}$ so that for any $1 \leq j_i \leq s_i$ and  $1\leq p^i_1<p^i_2< \cdots <p^i_{j_i} \leq k_i \leq l_i$, $i=1,2, \cdots ,r-1$, the intersection $T_{a_1, a_2} \cap \cdots \cap T_{a_1, a_k}$ has dimension at most $s-j$. For now, fix $1 \leq j_i \leq s_i$ and $1\leq p^i_1<p^i_2< \cdots <p^i_{j_i} \leq k_i \leq l_i$, $i=1,2, \cdots ,r-2$,      $1 \leq j_{r-1} \leq s_{r-1}$ and $1\leq p^{r-1}_1<p^{r-1}_2< \cdots <p^{r-1}_{j_{r-1}} < k_{r-1} \leq l_{r-1}$. By the induction hypothesis, that $T_{a_1, a_2} \cap \cdots \cap T_{a_1, a_{j-1}}$ has dimension at most $s-j+1$.

Split the variety $T_{a_1, a_2} \cap \cdots \cap T_{a_1, a_{j-1}}$ into irreducible components $W_1, \ldots, W_r$ and suppose that $W_a$ is a component of dimension $s-j+1 \geq 1$. By Lemma $2.2,  W_a$ has infinitely many points. Fix $d+1$ points $w_1, \ldots, w_{d+1}$ on $W_a$. For any $(s-1)$-variate polynomial $f $, write
$$
T_{a_1, f}=\left\{\left(x_1, \ldots, x_{s-1}\right):\left(f-f_{a_1}\right)\left(x_1, \ldots, x_{\prod_{i=1}^{r-1} s_i-1}\right)=0\right\} .
$$

By Lemma 2.3, we see that if $\operatorname{dim} W_a \cap T_{a_1, f}=\operatorname{dim} W_a$, then $T_{a_1, f}$ must contain all of $W_a$ and, in particular, each of $w_1, \ldots, w_{d+1}$. Therefore, for a random $(s-1)$-variate polynomial $f$ of degree $d$, the probability that $W_a \cap T_{a_1, f}$ does not have dimension at most $s-j$ is at most the probability that the polynomial $f  -f_{a_1}$ passes through all of $w_1, \ldots, w_{d+1}$, which, by Lemma $2.1$, is at most $q^{-(d+1)}$.

Since, by Lemma 2.4, the number of irreducible components of $T_{a_1, a_2} \cap \cdots \cap T_{a_1, a_{j-1}}$ is at most $d^{s-1}$, this implies that the probability $T_{a_1, a_2} \cap \cdots \cap T_{a_1, a_{j-1}} \cap T_{a_1, a_k}$ does not have dimension at most $s-j$ is at most $d^{s-1} q^{-(d+1)}$. By taking a union bound over the at most $\ell^{s-1}$ choices for $j$ and $a_1,a_2,\cdots,a_{j-1}$, so the probability there exists $1 \leq j_i \leq s_i$ and  $1\leq p^i_1<p^i_2< \cdots <p^i_{j_i} \leq k_i \leq l_i$, $i=1,2, \cdots ,r-1$ such that $T_{a_1, a_2} \cap \cdots \cap T_{a_1, a_{j-1}} \cap T_{a_1, a_{k}}$ does not have dimension at most $s-j$ is at most $\ell^{s-1} d^{s-1} q^{-(d+1)}<1$ for $q$ sufficiently large. Therefore, there exists an $(s-1)$-variate polynomial $f_{a_k}$ of degree at most $d$ such that $T_{a_1, a_2} \cap \cdots \cap T_{a_1, a_{j-1}} \cap T_{a_1, a_{k}}$ has dimension at most $s-j$  for any $1 \leq j_i \leq s_i$ and  $1\leq p^i_1<p^i_2< \cdots <p^i_{j_i} \leq k_i \leq l_i$, $i=1,2, \cdots ,r-1$.
\end{proof}

\end{document}